\documentclass{article}
\usepackage[utf8]{inputenc}
\usepackage{xcolor}
\usepackage{lscape}
\usepackage{array}
\usepackage{hyperref}

\usepackage[margin=1.5in]{geometry}
\usepackage{comment}
\usepackage[inline]{enumitem}
\usepackage{amsthm}

\usepackage{soul}

\usepackage[colorinlistoftodos]{todonotes}

\title{The Mathematics of Mathematics: Using Mathematics and Data Science to Analyze the Mathematical Sciences Community and Enhance Social Justice}

\author{Ron Buckmire \and Joseph E. Hibdon, Jr. \and Drew Lewis  \and Omayra Ortega \and Jos\'e L. Pab\'on \and Rachel Roca \and Andr\'es R. Vindas-Mel\'endez}

\begin{document}

\maketitle
%%%%%%%%%%%%%%%%%%%%%%%%%%%%%%%%%%%%

\begin{abstract}
We present and discuss a curated selection of recent literature related to the application of quantitative techniques, tools, and topics from mathematics and data science that have been used to analyze the mathematical sciences community.
We engage in this project with a focus on including research that highlights, documents, or quantifies (in)equities that exist in the mathematical sciences, specifically, and STEM (science, technology, engineering, and mathematics) more broadly.
We seek to enhance social justice in the mathematics and data science communities by providing numerous examples of the ways in which the mathematical sciences fails to meet standards of equity, equal opportunity and inclusion.
We introduce the term ``mathematics of Mathematics” for this project, explicitly building upon the growing, interdisciplinary field known as ``Science of Science” to interrogate, investigate, and identify the nature of the mathematical sciences itself.
We aim to promote, provide, and posit sources of productive collaborations and we invite interested researchers to contribute to this developing body of work.

\end{abstract}
%%%%%%%%%%%%%%%%%%%%%%%%%%%%%%%%%%%%

\section{Introduction}
\label{sec:intro}

\subsection{Motivation}
\label{sub:motivation}

Mathematics has long been used to solve problems, document patterns, describe phenomena, explain observations, and provide insight into how the world works. 
In the last century, with the advent of the computing age and its concomitant generation and storage of vast quantities of data, the ability to analyze data has become another important tool of mathematicians, engineers, and scientists.
 The context for what kinds of problems mathematics and data science can be used to solve are extremely varied, running the gamut from those that are theoretical with no currently foreseen applications to those that are immediately applicable to real-world situations. 
In this paper we are interested in how tools, topics, and techniques from the mathematical sciences can be and have been used to analyze the mathematical sciences itself. 
We introduce here the phrase ``the mathematics of Mathematics" (or \#MetaMath, for short) to refer to this growing body of work using mathematics to analyze Mathematics. 
Note that we are using the lower-case term ``mathematics'' in this phrase as shorthand for quantitative and qualitative techniques commonly associated with the mathematical sciences, and the upper-case term ``Mathematics" to mean the people, practices, and organizational structures associated with the profession of mathematics. In particular, we seek to highlight that mathematics is done by people who exist as part of a larger community, not in a vaccuum.  
We chose the term ``mathematics of Mathematics" specifically as an allusion to the growing field of ``science of Science"\footnotemark \cite{fortunato2018science,wang2021science,zeng2017science}, which shares many of our goals.

\footnotetext{We note that despite the similarity in nomenclature, neither of these terms have any association with the so-called ``science of math'' movement, a reactionary movement in K-12 education that rejects many equity-minded reforms that are based on mathematics education research. We further note that \#metamath has no connection to metamathematics which is a branch of mathematical logic.}

\subsection{Distinguishing  Mathematics of Mathematics (\#MetaMath) from Science of Science (SciSci) and Other Similar Disciplines}
%\cite{fortunato2018science,wang2021science,zeng2017science}
The science of science (SciSci) is a well-developed field of study in which the tools and techniques of scientific exploration are used to investigate and interrogate science itself. For example, Fortunato et al. \cite{fortunato2018science} say 
\begin{quote}
    SciSci places the practice of science itself under the microscope, leading to a quantitative understanding of the genesis of scientific discovery, creativity, and practice and developing tools and policies aimed at accelerating scientific progress.\\

    \noindent [...]\\

    \noindent SciSci relies on a broad collection of quantitative methods, from descriptive statistics and data visualization to advanced econometric methods, network science approaches, machine-learning algorithms, mathematical analysis, and computer simulation, including agent-based modeling.
\end{quote}
We intend the mathematics of Mathematics (\#MetaMath) to take advantage of the increased availability of data about different aspects of the Mathematics community, such as public information about research publications, citation patterns, faculty hiring, racial, ethnic,  and gender demographics, doctoral degree advisor-advisee relationships, research funding, and more, in order to explicitly advance the social justice goal of a more equitable, diverse, and inclusive Mathematics community, profession, and enterprise. 
(We note that just as the science of science includes qualitative techniques we intend for the mathematics of Mathematics to do so as well.)

\subsection{Context and Goals For This Paper}
This review paper is written in the context of recent activity in  how data science and mathematics are being used to enhance social justice in various new and exciting ways. 
For example, we envision \#MetaMath as a \emph{quantitative justice} project, where the term ``quantitative justice" is defined as ``the application of techniques, tools and topics from the quantitative sciences (e.g., mathematics, applied mathematics, data science, computer science, etc.) in subject domains that are derived from and/or typically associated with the social sciences (e.g., history, political science, law, economics, sociology etc.) with the explicit goal of promoting social justice” \cite{buckmireAAAS2023,buckmireAN2022qj,buckmireJMM2023metamath,buckmireCSE2023,buckmireJMM2023qj} .

Talks in the developing area of quantitative justice have occurred at multiple minisymposia \cite{SIAMAN22QJ} and special sessions \cite{SIAMJMM23} at several widely attended international STEM conferences in the last two years.
Examples include the 2022 annual meeting of the Society for Industrial and Applied Mathematics (SIAM) in Pittsburgh, PA \cite{buckmireAN2022qj}, the 2023 SIAM Computer Science and Engineering Conference in Amsterdam, the Netherlands \cite{buckmireCSE2023}, the 2023  annual meeting of the American Association for the Advancement of Science (AAAS) in Washington, DC \cite{buckmireAAAS2023}, and the 2023 Joint Mathematics Meetings in Boston, MA \cite{brooks2023branching,buckmireJMM2023qj,buckmireJMM2023metamath}.

In addition to being described as a quantitative justice project, \#MetaMath   also meets the definition of Data Science for Social Justice (DS4SJ)  provided in \cite{jones2023DS4SJ} as ``data scientific work (broadly construed) that actively
challenges systems of inequity and concretely supports
the liberation of oppressed and marginalized
communities."
Regardless of which term  one uses to describe the \#MetaMath project (quantitative justice, data science for social justice, justice data science, or something else), its goal is the same as this paper: to enhance social justice by using the mathematical sciences to investigate, interrogate, and identify inequities in the mathematical sciences itself.

In the main body of the paper, we provide a survey of selected recent work that uses tools, topics, and techniques from mathematics and data sciences to describe, document, and display inequities in the mathematics and data science communities. In the rest of the paper, we discuss some methodological considerations for quantitative scientists and mathematicians interested in working in the field, as well as some future directions for research.

\section{Survey of Existing \#MetaMath Research}
\label{sec:survey}

In this section we will present our commentary on  various examples of articles that have appeared in  the research literature over the last two decades that can be characterized with a \#MetaMath label, i.e. they use tools, topics, or techniques from mathematics and data science to analyze various aspects of discrimination, disparity, and difference in the mathematics community specifically, or the STEM academic community, generally.
We note that our presentation is not an exhaustive list or review of all such examples, since we are particularly interested in providing evidence that describes, displays, or documents inequity in the mathematics community.

If one believes the saying, as we do, that “Talent is evenly distributed, but opportunity is not,” then the research presented below bolsters that view. 
It shows that there exist hierarchical structures in mathematics, evidence of elitism, disparities in distribution of resources, over-representation of certain favored groups, under-representation of other disfavored groups, and a paucity of data on mathematics curriculum. 
We hope by highlighting \#metamath research that uses data to documents these deficiencies we will encourage others to enhance social justice in the mathematical sciences community by enacting positive change.

In the subsequent subsections we provide several examples of peer-reviewed  articles that use mathematics and data science to demonstrate these phenomena exist in STEM broadly and mathematics specifically.
We note that these phenomena are interrelated and interdependent  so it can be difficult to assign a particular paper or result to  only one of these overlapping areas.
For example, research by Liu et al \cite{liu2023citations} that shows that women and racial/ethnic minorities receive fewer citations to their research in mathematics and scientific journals demonstrates the lack of diversity in these fields, but also the existence of a hierarchy based on race, ethnicity and gender in academic publishing as well as be indicative of elitism in this area.
In cases like this the same research may be referred to in multiple parts of the paper.

\subsection{Diversity and Demographics of the Mathematics Community}

We believe it is self-evident that mathematicians, i.e. the members of the mathematics community in the United States (or however one defines the term \cite{buckmire2023definitions}) do not represent the full spectrum of the diversity and demographics of the general population. 
In other words, there are groups that are currently under-represented and historically marginalized in the mathematics community \cite{buckmire2021diversity}.

This under-representation of women, racial and ethnic minorities, people whose parents or grandparents didn't go to college (i.e., first-generation), people living with disabilities and members of the LGBTQ community has been documented and analyzed in many STEM fields, including in mathematics, and in different professional activities in academia.
The sectors of academia activity in which this under-representation manifests is varied and near-ubiquitous.
In other words, it is difficult to find an aspect of academia in which the diversity and demographics of those who participate in that aspect match the corresponding racial, ethnic and gender breakdowns of the general population of the United States.

Data on the demographics of the STEM community can be found in documents created by the National Center for Education Statistics in the U.S. Department of Education (NCES) \cite{nces2021digest} and the National Science Foundation's National Center for Science and Engineering Statistics (NCSES) \cite{ncses2023wmpd}. For math-specific data, the quinquennial survey conducted every five years since 1970 by the Conference Board of the Mathematical Sciences (CBMS) and the annual survey of earned doctorates by the American Mathematical Society (AMS) are key data sources\footnotetext{The CBMS survey can be found at \url{https://www.ams.org/profession/data/annual-survey/annual-survey} and the Annual Survey of Earned Doctorates can be found at \url{https://www.ams.org/profession/data/annual-survey/annual-survey}. Due to the COVID-19 pandemic, updates of these data sources has been significantly delayed.}

Buckmire \cite{buckmire2021diversity} notes that between 2013 to 2018, women made up approximately 42.6\% of recipients of bachelor's degrees in Mathematics. 
In the same time period, the racial and ethnic makeup of undergraduate math graduates ranged from 64.9\% White, 7.5\% Latino/Hispanic, and 5.0\% Black to 52.5\% White, 9.9\% Latino/Hispanic and 4.2\% Black.

Medina \cite{medina2004doctorate} notes that between 1993 and 2002, less than 5\% of those who earned doctoral degrees in mathematics were Black, Latinx, or Indigenous, despite those communities composing a quarter of the US population at that time.
Martin \cite{martin2014annotated} compiled a comprehensive list of 95 citations focusing on women in mathematics in his annotated bibliography, covering biases and inequities in primary school to tenure.
Cech \cite{Cechlgbtqstem2021} found that LGBTQ+ scientists and mathematicians were  more likely to report social marginalization, work place harassment, as well as limited career opportunities.

Topaz et al. \cite{topaz2016gender} analyzed 13,067 editorial positions on the boards of 435 mathematical science journals and found that roughly 8.9\% of editorships were held by women.  More recently, Liu et al. \cite{liu2023citations} expanded on this analysis in multiple dimensions. They looked at over one million papers between 2001 to 2020 in over 500 journals processed by almost 65,000 editors and found that ``non-white scientists appear on fewer editorial boards, spend more time under review, and receive fewer citations" than White scientists. 

In this subsection we have presented examples of data that documents how the demographics of the mathematics (and wider STEM) community do not reflect the diversity of the larger population from which it draws.
In the next subsection, we will present examples of research that addresses the similar but distinct problem that mathematics (and STEM) contain hierarchical structures that reflect and reproduce the lopsided demographic distributions previously discussed above.

\subsection{The Myth of Meritocracy and Existence of Hierarchy in Mathematics and Science}
%{\color{red}\cite{clauset2015systematic, wapman2022quantifying,fitzgerald_huang_leisman_topaz_2022,kaminski2012survival,mihaljevic2016effect,kawakatsu2021emergence,mihaljevic2020authorship,myers2011mathematical,morgan2022socioeconomic}}

Here we present examples of published research articles that can be used to prove the existence of hierarchical structures in mathematics (and STEM), thus refuting the pernicious claim that mathematics (or STEM) is a meritocracy.
Examples of hierarchical structures \cite{kawakatsu2021emergence} could include the existence of institutions which have greater prestige \cite{myers2011mathematical} or are more likely to have students  who are more likely to go on to obtain faculty positions at more prestigious institutions.

In the last few years, several researchers have taken advantage of the availability of vast quantities of data  on faculty positions at institutions of higher education in the United States to document the existence of hierarchies in faculty hiring networks in academia. 
Clauset et al. \cite{clauset2015systematic} demonstrated the existence of hierarchy in faculty hiring in a study involving 19,000 faculty members in a total of 461 departments in computer science, business, and history.  Wapman et al. \cite{wapman2022quantifying} expanded this analysis to cover 295,089 faculty in 10,612 departments at 368 Ph.D.-granting institutions and all academic disciplines for the 2011-2020 decade. FitzGerald et al. \cite{fitzgerald_huang_leisman_topaz_2022} built upon this research by using data from the Mathematics Genealogy Project (MGP) to restrict their analysis to mathematics faculty. 
They analyzed  120,000  records from 150 institutions over seven decades (from 1950 to 2019). 
The results of this research demonstrates that who becomes faculty at  Ph.D.-granting institutions is not an unbiased process or reflective of a meritocratic ideal.
Instead, analysis of the availability data shows that there exist certain institutions that are more likely than others to have their Ph.D. graduates  become faculty members at other Ph.D.-granting institutions.
In other words, the data demonstrates that hierarchies exist in faculty hiring networks.
%Mucha et al. \cite{myers2011mathematical} also analyzed data from the MGP to 

One of the most salient and persistent hierarchical structures in society that is omnipresent in mathematics and the wider STEM community is gender.
Researchers have analyzed data describing many different aspects of academic activity and demonstrated the many ways gender can negatively mediate opportunity for advancement, participation, and achievement in science and mathematics \cite{brisbin2015women,kaminski2012survival,martin2014annotated,mihaljevic2020authorship,mihaljevic2016effect,rissler2020gender,schmaling2023gender,topaz2016gender,way2016gender}.

Another significant hierarchical structure in society is class. 
Morgan et al \cite{morgan2022socioeconomic} analyzed data from a survey with responses from more than 7,000 faculty in eight STEM disciplines to demonstrate that faculty are 25 times more likely to have a parent with a Ph.D. and that at ``elite" institutions this rate is doubled. 
This data analysis demonstrates the salience of socioeconomic status on the composition of the professoriate.

In this subsection, we have provided articles that confirm the existence of hierarchical structures in mathematics and science.
Our argument is that these research results using data demonstrate the idea that mathematics (and STEM, more generally) are egalitarian and meritocratic is a myth that can be and has been disproved, debunked, and discredited.

\subsection{Evidence of Elitism in Mathematics} 

%Removing the citations from this subsection title.
%\cite{chang2021elitism,schlenker2020prestige,kawakatsu2021emergence,mihaljevic2016effect,topaz2016gender,leslie2015expectations,storage2016frequency,way2016gender,brisbin2015women,mihaljevic2020authorship}

Our survey of existing results finds extensive evidence, gathered via mathematical tools and techniques, of widespread elitism and exclusion of minorities and under-represented populations in the mathematics community \cite{chang2021elitism, kawakatsu2021emergence,mihaljevic2016effect,schlenker2020prestige,topaz2016gender}. 
For example, Topaz et al. found in 2016 that in academia, women comprised only 15\% of tenure-stream faculty positions in doctoral-granting mathematical sciences departments in the United States. 
Editorial appointments in the journals studied in their work found even more paltry representation for women; only 8.9\% of editorships within the journals reviewed belonged to women \cite{topaz2020institute,topaz2016gender}. 

Another way elitism is perpetuated in the mathematics community is via the selection of a self-perpetuating cadre of ``elite" personnel for prestigious prizes. 
%Our global and local communities incentivize for or deter elitism in our societies via the programs, protocols and procedures they enact to command and control the socio-educational hierarchies they employ. 
The International Mathematical Union's selections for the Fields Medal, (widely considered the most elite prize in Mathematics) for example, have shown an alarming trend of rewarding only already-elite mathematics researchers in an apparent negation of one of the main goals of its conferral, which was to shine light on under-represented members of the mathematics community \cite{chang2021elitism}. 
Other scientists have found widespread elitism perpetuates as a vicious circle of like begetting like \cite{schlenker2020prestige, storage2016frequency}. 
%Prize and funding initiatives by existing entities in the mathematical community have failed to succeed in both raising the academic prestige status of the funded research and fostering a larger share of under-represented mathematicians in those fields.

Schlenker \cite{schlenker2020prestige} notes that fields with applications to the social or physical sciences such as numerical analysis, mathematical modeling or statistics seem to be viewed as having low status, and this lack of prestige accompanies the low representation of researchers in these fields among elite prizewinners.  
However, other fields such as Algebraic Geometry have no shortage of practitioners and medalists, in spite of having little direct application to other disciplines, once again displaying the self-serving cycle of elitism within the mathematics community.

Other researchers have also used other mathematical tools to identify inequity in mathematics in other aspects of academia. 
Leslie et al. \cite{leslie2015expectations} and  Storage et al. \cite{storage2016frequency} found widespread under-representation of women and African Americans in their research which analyzed over 14 million records from the popular website RateMyProfessor. 
Their findings included detailing that in fields where brilliance and innate intellectual talent were treasured, African Americans and women were even less well represented. 
Brisbin and Whitcher \cite{brisbin2015women}  concur with the above research by finding that women are under-represented in subfields of mathematics that are viewed as having more prestige by analyzing almost a million papers in the mathematical sciences uploaded to the arXiv, a popular website for the dissemination of research results.  
Way et al. \cite{way2016gender}  studied comprehensive data on both hiring outcomes and scholarly productivity for 2659 tenure-track faculty across 205 Ph.D.-granting departments in North America and found deep under-representation of women in these groups.  
Curiously, Way et al. also found pervasive elitism including discrimination against females, as well as some other less prominent examples of biases against women.  
One of the salient trends that these researchers found was that top-ranked institutions had tendencies to hire women at higher rates than their mathematical data science enabled models would expect. 
In this sense, higher ranked academic institutions did a ``better job'' of being less discriminating against women. 
This was also the case at the lower end of the spectrum of lower ranked academic institutions; thus the majority of hiring and appointing of women in academia take place in the middle ranked academic institutions.

The research highlighted in this subsection has provided examples of the pervasive elitism across academia at large and the mathematics community in particular.

\subsection{Disparities in Federal Funding for Science and Mathematics} %\cite{chenNSFdisparities,ginther2011race,rissler2020gender,taffe2021racial,schmaling2023gender}
Another example of \#MetaMath research involves analyzing the public data of who receives funding from public and private sources in order to detect and document disparities in support for mathematics and science in the United States.
In this subsection we provide examples of research that provides evidence of these discriminatory results.
It is well documented that there are racial disparities in funding from the National Institutes of  Health (NIH)  \cite{chenNSFdisparities,ginther2011race,rissler2020gender,schmaling2023gender,taffe2021racial}. 
The 2011 study of NIH R01 funding by Ginther et al. \cite{ginther2011race} hypothesized that scientists of different races but with similar research records and institutional affiliations would also have a similar likelihood of being awarded research grants. 
However, the opposite was discovered.
They found that Black and Asian investigators are less likely to be awarded an R01 on the first or second attempt, Blacks and Hispanics are less likely to resubmit a revised application, and Black investigators that do resubmit have to do so more often to receive an award. 
Even after controlling for the applicant’s educational background, country of origin, training, previous research awards, publication record, and employer characteristics, the authors discovered that black applicants remain 10 percentage points less likely than whites to be awarded NIH research funding. 

Funding  from the National Science Foundation (NSF) have not been as extensively studied as the NIH. 
As the main source for research funding in the mathematical sciences, the funding priorities selected by the NSF will greatly influence how the fields in the mathematical sciences will grow, behave, and look. 
Both the NIH and NSF are federal funding agencies with similar missions but different foci. 
As such, we would not expect to see  trends in NSF funding differ significantly from the NIH. 
Indeed, Chen et al. \cite{chenNSFdisparities} recently examined NSF funding data from 1996 to 2019 and found that proposals by white PIs were consistently funded at rates higher than the overall rate, with an average relative funding rate of +8.5\% from 1999 to 2019. Additionally, the relative funding rate for proposals by white PIs increased steadily during this period, from +2.8\% in 1999 to +14.3\% in 2019. 
Proposals by Asian, Black/African-American, and Native Hawaiian/Pacific Islander PIs, were consistently funded below the overall rate, with average relative funding rates of –21.2\%, –8.1\%, and –11.3\%, respectively.

Additionally, these findings revealing race-based  and gender-based funding disparities are not confined to federal agencies. 
The Wellcome Trust is a private global charitable foundation that conducted a similar study on their own award practices in 2020 \cite{Wellcome}. 
The Wellcome study stated that, ``As a funder, we know from our data that BAME\footnote{Black, Asian, and minority ethnic, also known as “BAME,” is an umbrella term, common in the United Kingdom, used to describe non-white ethnicities.}, and especially Black, applicants are less likely to be awarded Wellcome research grants in the UK than White applicants.'' 
At a time when the mathematics community has been doubling down on its efforts to diversify the profession, the type of groundbreaking data analysis presented in this subsection is absolutely necessary to motivate action and meaningful change to eliminate racial and gender disparities in funding  in the mathematical sciences community.

\subsection{Discourse in the Mathematics Community}
\label{sub:discourse}
Just like SciSci, \#MetaMath also includes the use of qualitative techniques.
In this subsection we provide some examples of the ways that discourse in and about the mathematics community can be analyzed, primarily through mixed or qualitative methods, in order to raise awareness about social justice issues in the mathematical sciences.

The mathematics community shares its opinions and responses to given situations in a variety of public ways, often leaving behind artifacts that can be studied. 
Examples of such dissemination of opinions and responses include the Inclusion/Exclusion blog \cite{inclusionexclusionblog}, a justice and mathematics weblog that is a continuation of the AMS blog archived in 2021; an article in \textit{Scientific American} highlighting the debate surrounding the founding of the Association for Mathematical Research \cite{Crowell2022AMR}; and the statement by the Mathematical Association of America (MAA)  on the location of Mathfest 2023 \cite{Mathfest2023}, in which they prioritize financial concerns over the safety of the LGBTQ+ community at the meeting site in a state that has enacted discriminatory legislation targeting LGBTQ+ people. The artifacts left behind in these cases include the blog, the article and the statement as well as any public replies to these on social media, statements from organizations, letters to the editor, and more. 
By examining the conversation through these artifacts, we can also better understand and quantify the divisions and counter-movements within the mathematics community. It is worth highlighting that this type of work often requires mixed-methodology; strictly quantitative analysis cannot fully capture nuance and ``whys'' behind the discourse.

One of the first questions we might ask ourselves in this work is \textit{who} is part of the conversation, which can be answered quantitatively. The work of Topaz et al. \cite{topaz2020comparing} compares the demographics of the hundreds of signatories on three public letters responding to an essay published in \textit{Notices of the American Mathematical Society} that opposed the use of diversity statements in academic hiring, likening them to ``McCarthyism.'' The results showed that signatories of Letter A, which highlights diversity and social justice, were inferred to be more diverse in terms of gender, under-represented ethnic groups, professional security, and employment at academic institutions. Letters B and C, which did not mention diversity and spoke against diversity statements, respectively, were instead signed by those who were inferred to be majority tenured white men at research focused universities. The results presented in this paper are consistent with theories of power and positionality found in social sciences \cite{merriam2001power, misawa2015cuts}. 

From a more qualitative perspective, Burrill et al. \cite{burrill2023notices} identified broad themes through recordings from a math education forum, consolidating and sharing common ground across diverse experts. Their work highlights the importance of communication and articulation among stakeholders, as well as the need to support students and teachers in curricular design. Thematic qualitative analysis can also be executed in a more fine-grained manner, such as through the examination of Twitter posts containing the \#disruptJMM hashtag by Roca et al. \cite{roca2023disruptjmm}. This paper examines the math social movement born from Piper H's post on the Inclusion/Exclusion Blog calling for more commonplace discussion on inclusion and equity within the Joint Mathematics Meeting (JMM) under the \mbox{\#disruptJMM} banner. Those who used the hashtag and challenged others to get involved built a community that highlighted the need to rehumanize mathematics, make visible power and privilege in the mathematics community, and elevate the stories of under-represented folks at the conference.

The examples of the ways discourse in the mathematics community can be studied given in this subsection allows us to identify salient topics of conversation, understand who is included within or excluded from these discussions, and provide context to situations affecting the heterogeneous mathematics community in disproportionate way.

\subsection{Quantifying Diversity in Mathematical Content  and Curricular Choices}
\label{sub:curriculum}

This subsection is motivated by two questions.  
First, \emph{what courses are mathematics departments currently implementing in their undergraduate programs}?
Second, \emph{how are mathematics departments using data to make curricular choices?}
The reality is that there is limited data available to properly answer these questions. 

Every five years since 1970, CBMS sponsors a comprehensive national survey of undergraduate mathematical and statistical sciences in the United States' four-year and two-year universities and colleges. 
The national data collected includes statistics on enrollment, curriculum, degrees awarded, course availability, faculty demographics, and special one-time topics (these topics vary) \cite{CBMS2015}. 
The most recent CBMS survey, whose complete results are published, is from 2015. (The release of data from the 2020 survey has been delayed due to the COVID-19 pandemic.)
This survey showed that among the 28 upper-level math courses listed, the courses that were most offered (at least once) in the academic years 2014-2016 and 2015-2016 were Modern Algebra I (78\% of all math departments),  Geometry (71\%), Advanced Calculus/Real Analysis I (72\%), Math Senior Seminar/Independent study (66\%), History of Mathematics (58\%), and Introduction to Proofs (56\%).

Although there are guidelines written by leading mathematical societies on what the mathematics undergraduate curriculum should consist of (see for instance \cite{GAISE,CUPM,DataUndergrad,GoldKeithMarion,GuidelinesUM, SelfStudy,  SaxeBraddy}), the CBMS survey demonstrates that there is not a consistent vision for the mathematics curriculum among the institutional respondents. 
Additionally, we must take into consideration the limitations of such surveys/data, such as survey response rates. 
For example, among the four-year mathematics respondents of the 2015 CBMS survey, only 5 out of the 23 California State University System (CSU) campuses and 4 out of the 9 University of California (UC) System undergraduate campuses participated. 
Combined, both the CSU and UC systems would create the largest public institution of higher education in the United States, thus there is a large number of students and departments not voiced in the survey results.

This leads us to the following question: \emph{why should there be data on mathematics course offerings at different institutions?}
Answering this question provides an avenue for at least two data science projects that can benefit the broader mathematical community. 
One such project would be to actually collect the data on course offerings at all institutions with an undergraduate mathematics major. 
This may involve applying data science techniques, such as data/web-scraping, information retrieval, and/or data clustering.
Having data on mathematics course offerings at different institutions can be used to influence funding, administrative, and educational practices. 
As far as we know, there is also no known information or literature on how mathematics departments are using quantitative data as an input when a department or faculty chooses math content and curricula.
Hence, a second project is to pursue a study towards understanding the decision-making processes of mathematics departments in regards to the aforementioned setting.
Such a study can have implications for mathematical practice and education and can help interested parties to identify related important and timely matters, such as textbook selection for courses. 

A deeper dive into the content of mathematics courses would inform the question, \emph{what is  mathematics?} This question is explored through the Rehumanizing Mathematics framework \cite{Gutierrez}, developed by Dr. Rochelle Gutierrez, and can have implications for future data science inquiry. 
In this framework, there are eight dimensions; participation/positioning; cultures/histories; windows/mirrors; living practice; broadening maths; creation; body/emotions; and ownership. 
These dimensions can be implemented in the process of rehumanizing the mathematics classroom. 
This framework has also been used in conjunction with lesson study to rehumanize the content of mathematics courses by bringing back diverse mathematics history that is often left out of mathematics courses, including and incorporating students' existing knowledge and experience into the classroom, and by considering alternative formats for conducting and communicating mathematics (for example). 
Thus far, no one has quantified and categorized the content offered in mathematics courses. 
This inquiry could be further expanded to identify geographical trends in the types of mathematics content offered.

Decisions about who gets admitted to college, particularly to elite or selective institutions, are the result of a combination of policy and practice, and math expectations play a role \cite{AndersonBurdman}. 
Strongly-held beliefs about calculus as a sign of rigor play a consequential role in college admissions.
This can also translate to upper-division mathematics courses and to graduate admissions.
Data, such as that presented in the CBMS survey or potential results from the data science projects proposed above, should be carefully used and analyzed.  
For example, graduate admission committees when reviewing a student's application should take into consideration questions such as: Did the student have opportunities to enroll in graduate courses as an undergraduate student? Did the student study at small college or large university? Did the student participate in research or an internship? Did the student have the opportunity to pursue an undergraduate thesis? 
%These results also come with their limitations and can lead to uncovering underlying educational and socioeconomic questions, or in some cases, may overlook or not take into account those situations. 
%Moreover, seeing the higher percentages of courses being offered in mathematics departments, such as those presented in the 2015 CBMS, may lead people to view Modern Algebra or Real Analysis as a sign of rigor, which can also play a consequential role in situations such as graduate admissions. 
The call in this subsection for the ability to quantify the diversity of mathematical content and curricular choices at undergraduate institutions could lead to more nuanced admissions decisions at graduate schools and more effective undergraduate mathematics curricula with a possibility of a future mathematics community that has been positively impacted by the existence and analysis of such data.

\section{Methodological Considerations}
\label{sec:method}

As members of the mathematical community embark on  work in \#MetaMath, we caution them to proceed with humility. Much of this work, while quantitative in nature, is epistemologically closer to the social sciences than to proving theorems or writing computer code. 
As such, we point out some common pitfalls we urge people considering this line of research to consider.

One common research approach is to analyze an existing public data set (\cite{brisbin2015women,brooks2023branching,clauset2015systematic,fitzgerald_huang_leisman_topaz_2022,reys2022some,schlenker2020prestige,wapman2022quantifying}).  
While this can provide important insights, this can be limiting in several ways.  
First, this approach inverts the typical research process, in which a research question is formulated, and then data obtained; here, data is obtained, and research questions are formed (implicitly or explicitly). 
In particular, this can limit the kinds of questions that can be asked, and as a result limit the kinds of conclusions that can be made. 
Second, data collection is not a neutral endeavor.  
Many choices are made in the production of a data set, and these choices can have profound implications for the kinds of questions that can be answered and the kinds of conclusions that can be drawn. For example, \cite{fitzgerald_huang_leisman_topaz_2022} makes use of the data contained in the popular Mathematics Genealogy Project. However, as the authors note, this has the effect of excluding from their dataset all mathematicians who have not advised a Ph.D. student.

Many of the studies given in previous subsections above concern questions about racial and gender representation in the mathematics community. However, the needed demographic information to conduct this data analysis is often not present.  This leaves researchers to collect new data, or to somehow infer demographic information from other information present in the data set. One approach is to try to match the names of individuals to other public information, such as a personal website, as in \cite{topaz2020comparing,topaz2016gender}. This is, of course, somewhat labor intensive.  An automated approach is to use software tools that attempt to infer an individual's gender based on their name (and sometimes other information), as in \cite{fitzgerald_huang_leisman_topaz_2022}.  This practice of researchers ascribing gender or race to an individual is certain to contain errors to some degree, and should be used with caution---researchers should carefully consider whether the benefits of these techniques truly outweigh the harms caused by incorrectly ascribing identities to people. 
Further, many studies including gender as a factor completely disregard nonbinary individuals, either marking them as ``other'' or discarding them from the data set.

Mathematicians undertaking this work should also be wary of an epistemological bias.  
Our training understandably predisposes us to use quantitative techniques.  
However, many questions of interest to the mathematics community are better suited to mixed or qualitative methodologies, and in fact may be impossible to answer with only quantitative methods. 
For example, we above highlighted several papers that illuminate demographic disparities in various aspects of the mathematical community.  
The next important step is to then interrogate why these disparities exist and persist, which is much better suited to qualitative methodologies. As such, we encourage those engaging in this type of work to foster collaborations with experienced qualitative researchers.  
The authors' personal experience is that the mathematics education research community (in addition to other social science researchers) is a natural source of collaborators that can help push this work forward.

\section{Future Directions}
It should be clear to the reader that the \#MetaMath examples given above are not an exhaustive list of the ways that mathematics and data science can be used to examine the mathematical sciences itself. 
There are many different avenues that future work could build upon the information we have provided in this article.
First, a review of the points raised in Section~\ref{sec:method} on Methodological Considerations could lead to a number of ways to improve upon or modify and replicate some of the results discussed in the survey of the literature that we have provided.
Second, examining the number of references where the context for the research is STEM and not exclusively the mathematical sciences could provide opportunities to repeat STEM-specific research in Math-exclusive contexts.
Third, inspired by the subsection discussing qualitative analysis of discourse in the mathematics community, future researchers could try to expand the use of qualitative techniques  to analyze the mathematical sciences.
Fourth, there are numerous suggestions for data-enabled projects involving curricular data in the subsection on mathematical content and curricular choices that could benefit the wider mathematics community.
Lastly, as more and more data about the mathematical sciences community becomes available, either updated versions of existing data or data about new topics or in areas where data is currently incomplete or inaccessible, these sources will provide new opportunities for ways to provide insight into the nature of the mathematics community through the application of techniques, tools, and topics from mathematics and data science.

We conclude by inviting interested researchers to join us in this project. The work of promoting social justice and improving equity within the field of mathematics is and will remain an ongoing process.

\section*{Author Contributions} All the listed co-authors contributed equally to the creation of this article; the order of attribution is alphabetical, as is customary in mathematics and is not intended to demonstrate any distinction in credit or effort.

\section*{Acknowledgments}
This material is based upon work supported by the National Science Foundation under Grant No. DMS-1929284 while some of the authors were in residence at the Institute for Computational and Experimental Research in Mathematics in Providence, RI, during the \href{https://icerm.brown.edu/programs/ep-22-dssj/}{Data Science and Social Justice: Networks, Policy, and Education program}.
The authors thank ICERM for providing the space for this  work. 

The authors would like to acknowledge Katherine (Katie) Kinnaird, who contributed to early discussions.
 
 Buckmire acknowledges that the Data Science and Social Justice program at ICERM, under lead organizer Carrie Diaz Eaton,  provided key organizational and logistical resources that facilitated the relationships and connections that resulted in this paper. 

Hibdon was supported, in part, by the National Institutes of Health’s National Cancer Institute, Grant Numbers U54CA202995, U54CA202997, and U54CA203000. The content is solely the responsibility of the authors and does not necessarily represent the official views of the National Institutes of Health.

Pab\'on thanks Nicholas Dubicki for insightful discussions. Pab\'on is partially supported by the National Science Foundation under Awards DMS-2108839 and DMS-1450182 .

Vindas-Mel\'endez thanks Benjamin Braun for helpful discussion. 
Vindas-Mel\'endez is partially supported by the National Science Foundation under Award DMS-2102921.

%%%%%%%%%%%%%%%%%%%%%%%%%%%%%%%%%%%%

\bibliographystyle{plain}
\bibliography{bibliography}

\begin{thebibliography}{10}

\bibitem{AndersonBurdman}
Veronica Anderson and Pamela Burdman.
\newblock A new calculus for college admissions: How policy, practice, and
  perceptions of high school math education limit equitable access to college,
  2022.
\newblock
  https://justequations.org/resource/a-new-calculus-for-college-admissions-how-policy-practice-and-perceptions-of-high-school-math-education-limit-equitable-access-to-college.

\bibitem{CBMS2015}
Richelle~M. Blair, Ellen~E. Kirkman, and James~W. Maxwell.
\newblock {\em Statistical abstract undergraduate programs in the Mathematical
  Sciences in the United States: 2015 CBMS survey}.
\newblock American Mathematical Society, 2018.

\bibitem{brisbin2015women}
Abra Brisbin and Ursula Whitcher.
\newblock Women's representation in mathematics subfields: Evidence from the
  arxiv.
\newblock {\em arXiv preprint arXiv:1509.07824}, 2015.

\bibitem{brooks2023branching}
Heather~Z Brooks, Philip Chodrow, Harlin Lee, and Juan~G Restrepo.
\newblock Branching process model for female gender representation in the
  mathematics genealogy project.
\newblock In {\em 2023 Joint Mathematics Meetings (JMM 2023)}. AMS, 2023.

\bibitem{buckmire2021diversity}
Ron Buckmire.
\newblock “{W}ho {D}oes the {M}ath?”: On the diversity and demographics of
  the mathematics community in the {USA}.
\newblock In {\em Improving {A}pplied {M}athematics {E}ducation}, pages 1--12.
  Springer, 2021.

\bibitem{buckmireAN2022qj}
Ron Buckmire.
\newblock An introduction to quantitative justice.
\newblock In {\em 2022 Annual Meeting of the Society for Industrial and Applied
  Mathematics (SIAM AN22)}. SIAM, 2022.

\bibitem{buckmireJMM2023qj}
Ron Buckmire.
\newblock An introduction to quantitative justice: The application of
  quantitative techniques from data science and mathematics to promote social
  justice.
\newblock In {\em 2023 Joint Mathematics Meetings (JMM 2023)}. AMS, 2023.

\bibitem{buckmireAAAS2023}
Ron Buckmire.
\newblock \#metamath: Using big data to promote social justice in the
  mathematics community.
\newblock In {\em Annual Meeting of the American Association for the
  Advancement of Science (AAAS23)}. AAAS, 2023.

\bibitem{buckmireCSE2023}
Ron Buckmire.
\newblock A quantitative justice primer: The \#metamath project and other
  examples.
\newblock In {\em SIAM Conference on Computer Science and Engineering (SIAM
  CSE23)}. SIAM, 2023.

\bibitem{buckmireJMM2023metamath}
Ron Buckmire, Carrie~Diaz Eaton, Joseph~Edward Hibdon, Katherine~M Kinnaird,
  Drew Lewis, Jessica~M Libertini, Omayra Ortega, Rachel Roca, and Andres~R
  Vindas-Melendez.
\newblock The\# metamath project: Analyzing the mathematics community using
  mathematics and data science to promote social justice.
\newblock In {\em 2023 Joint Mathematics Meetings (JMM 2023)}. AMS, 2023.

\bibitem{buckmire2023definitions}
Ron Buckmire, Carrie~Diaz Eaton, Joseph~E Hibdon~Jr, Katherine~M Kinnaird, Drew
  Lewis, Jessica Libertini, Omayra Ortega, Rachel Roca, and Andr{\'e}s~R
  Vindas-Mel{\'e}ndez.
\newblock On definitions of ``mathematician''.
\newblock {\em arXiv preprint arXiv:2302.07432}, 2023.

\bibitem{SIAMAN22QJ}
Ron Buckmire, Suzanne Lenhart, and Suzanne Sindi.
\newblock {SIAM} quantitative justice: Intersections between quantitative
  sciences and social justice, 2023.
\newblock SIAM Annual Meeting, Pittsburgh, PA.

\bibitem{SIAMJMM23}
Ron Buckmire, Omayra Ortega, and Carrie Diaz~Eaton.
\newblock {SIAM} minisymposium on quantitative justice (a {NAM-SIAM} joint
  session), 2023.
\newblock Joint Mathematics Meetings, Boston, MA.

\bibitem{burrill2023notices}
Gail Burrill, Henry Cohn, Yvonne Lai, Dev~P. Sinha, Ji~Y. Son, and Katherine~F.
  Stevenson.
\newblock Listening for {C}ommon {G}round in {H}igh {S}chool and {E}arly
  {C}ollegiate {M}athematicss.
\newblock {\em Notices of the American Mathematical Society}, 60(5):798--805,
  2023.

\bibitem{Cechlgbtqstem2021}
E.~A. Cech and T.~J. Waidzunas.
\newblock Systemic inequalities for lgbtq professionals in stem.
\newblock {\em Science Advances}, 7(3):eabe0933, 2021.

\bibitem{chang2021elitism}
Ho-Chun~Herbert Chang and Feng Fu.
\newblock Elitism in mathematics and inequality.
\newblock {\em Humanities and Social Sciences Communications}, 8(1):1--8, 2021.

\bibitem{chenNSFdisparities}
Christine~Yifeng Chen, Sara~S Kahanamoku, Aradhna Tripati, Rosanna~A Alegado,
  Vernon~R Morris, Karen Andrade, and Justin Hosbey.
\newblock Meta-{R}esearch: {S}ystemic racial disparities in funding rates at
  the {N}ational {S}cience {F}oundation.
\newblock {\em eLife}, 11:e83071, 2022.

\bibitem{clauset2015systematic}
Aaron Clauset, Samuel Arbesman, and Daniel~B Larremore.
\newblock Systematic inequality and hierarchy in faculty hiring networks.
\newblock {\em Science advances}, 1(1):e1400005, 2015.

\bibitem{GAISE}
ASA~Revision Commitee.
\newblock Guidelines for assessment and instruction in statistics education
  ({GAISE}) in statistics education college report, 2015.
\newblock
  {https://www.amstat.org/docs/default-source/amstat-documents/gaisecollege\_full.pdf}.

\bibitem{Crowell2022AMR}
Rachel Crowell.
\newblock New math research group reflects a schism in the field.
\newblock {\em Scientific American}, 2022.
\newblock
  https://www.scientificamerican.com/article/new-math-research-group-reflects-a-schism-in-the-field.

\bibitem{nces2021digest}
C~De~Brey, TD~Snyder, A~Zhang, and SA~Dillow.
\newblock Digest of education statistics.
\newblock {\em National Center for Education Statistics, Institute of Education
  Sciences, US Department of Education. https://nces.ed.gov/programs/digest},
  2021.

\bibitem{DataUndergrad}
Richard~D. De~Veaux, Mahesh Agarwal, Maia Averett, Benjamin~S. Baumer, Andrew
  Bray, Thomas~C. Bressoud, Lance Bryant, Lei~Z. Cheng, Amanda Francis, Robert
  Gould, Albert~Y. Kim, Matt Kretchmar, Qin Lu, Ann Moskol, Deborah Nolan,
  Roberto Pelayo, Sean Raleigh, Ricky~J. Sethi, Mutiara Sondjaja, Neelesh
  Tiruviluamala, Paul~X. Uhlig, Talitha~M. Washington, Curtis~L. Wesley, David
  White, and Ping Ye.
\newblock Curriculum {G}uidelines for {U}ndergraduate {P}rograms in {D}ata
  {S}cience.
\newblock {\em Annual Review of Statistics and Its Application}, 4(1):15--30,
  2017.

\bibitem{fitzgerald_huang_leisman_topaz_2022}
Cody FitzGerald, Yitong Huang, Katelyn~Plaisier Leisman, and Chad~M. Topaz.
\newblock Temporal dynamics of faculty hiring in mathematics.
\newblock {\em Humanities and Social Sciences Communications}, 10(1):247, 2023.

\bibitem{ncses2023wmpd}
National~Center for Science and Engineering~Statistics (NCSES).
\newblock {\em Diversity and STEM: Women, Minorities, and Persons with
  Disabilities 2023 (Special Report NSF 23-315)}.
\newblock National Science Foundation, 2023.

\bibitem{fortunato2018science}
Santo Fortunato, Carl~T Bergstrom, Katy B{\"o}rner, James~A Evans, Dirk
  Helbing, Sta{\v{s}}a Milojevi{\'c}, Alexander~M Petersen, Filippo Radicchi,
  Roberta Sinatra, Brian Uzzi, et~al.
\newblock Science of {S}cience.
\newblock {\em Science}, 359(6379):eaao0185, 2018.

\bibitem{ginther2011race}
Donna~K Ginther, Walter~T Schaffer, Joshua Schnell, Beth Masimore, Faye Liu,
  Laurel~L Haak, and Raynard Kington.
\newblock Race, ethnicity, and {NIH} research awards.
\newblock {\em Science}, 333(6045):1015--1019, 2011.

\bibitem{GoldKeithMarion}
Bonnie Gold, Sandra~Z. Keith, and William~A. Marion, editors.
\newblock {\em Assessment practices in undergraduate mathematics}.
\newblock Notes 249. Mathematical Association of America, 1999.

\bibitem{Gutierrez}
R.~Gutiérrez.
\newblock {\em In Rehumanizing mathematics for Black, Indigenous, and Latinx
  students}, chapter Introduction: The need to rehumanize mathematics.
\newblock National Council of Teachers of Mathematics, 2018.

\bibitem{inclusionexclusionblog}
Inclusion/exclusion, a justice and math weblog.
\newblock https://inclusionexclusion.org/.

\bibitem{jones2023DS4SJ}
Quindel Jones, Andr{\'e}s R~Vindas Mel{\'e}ndez, Ariana Mendible, Manuchehr
  Aminian, Heather~Z Brooks, Nathan Alexander, Carrie~Diaz Eaton, and Philip
  Chodrow.
\newblock Data science and social justice in the mathematics community.
\newblock {\em arXiv preprint arXiv:2303.09282}, 2023.

\bibitem{kaminski2012survival}
Deborah Kaminski and Cheryl Geisler.
\newblock Survival analysis of faculty retention in science and engineering by
  gender.
\newblock {\em Science}, 335(6070):864--866, 2012.

\bibitem{kawakatsu2021emergence}
Mari Kawakatsu, Philip~S Chodrow, Nicole Eikmeier, and Daniel~B Larremore.
\newblock Emergence of hierarchy in networked endorsement dynamics.
\newblock {\em Proceedings of the National Academy of Sciences},
  118(16):e2015188118, 2021.

\bibitem{leslie2015expectations}
Sarah-Jane Leslie, Andrei Cimpian, Meredith Meyer, and Edward Freeland.
\newblock Expectations of brilliance underlie gender distributions across
  academic disciplines.
\newblock {\em Science}, 347(6219):262--265, 2015.

\bibitem{liu2023citations}
Fengyuan Liu, Talal Rahwan, and Bedoor AlShebli.
\newblock Non-white scientists appear on fewer editorial boards, spend more
  time under review, and receive fewer citations.
\newblock {\em Proceedings of the National Academy of Sciences},
  120(13):e2215324120, 2023.

\bibitem{martin2014annotated}
Greg Martin.
\newblock An annotated bibliography of work related to gender in science.
\newblock {\em arXiv preprint arXiv:1412.4104}, 2014.

\bibitem{medina2004doctorate}
Herbert~A Medina.
\newblock Doctorate degrees in mathematics earned by {B}lacks,
  {H}ispanics/{L}atinos, and {N}ative {A}mericans: a look at the numbers.
\newblock {\em Notices of the American Mathematical Society}, 51(7):772--775,
  2004.

\bibitem{merriam2001power}
Sharan~B Merriam, Juanita Johnson-Bailey, Ming-Yeh Lee, Youngwha Kee, Gabo
  Ntseane, and Mazanah Muhamad.
\newblock Power and positionality: Negotiating insider/outsider status within
  and across cultures.
\newblock {\em International journal of lifelong education}, 20(5):405--416,
  2001.

\bibitem{mihaljevic2020authorship}
Helena Mihaljevi{\'c} and Luc{\'\i}a Santamar{\'\i}a.
\newblock Authorship in top-ranked mathematical and physical journals: Role of
  gender on self-perceptions and bibliographic evidence.
\newblock {\em Quantitative Science Studies}, 1(4):1468--1492, 2020.

\bibitem{mihaljevic2016effect}
Helena Mihaljevi{\'c}-Brandt, Luc{\'\i}a Santamar{\'\i}a, and Marco Tullney.
\newblock The effect of gender in the publication patterns in mathematics.
\newblock {\em PLoS One}, 11(10):e0165367, 2016.

\bibitem{misawa2015cuts}
Mitsunori Misawa.
\newblock Cuts and bruises caused by arrows, sticks, and stones in academia:
  Theorizing three types of racist and homophobic bullying in adult and higher
  education.
\newblock {\em Adult Learning}, 26(1):6--13, 2015.

\bibitem{morgan2022socioeconomic}
Allison~C Morgan, Nicholas LaBerge, Daniel~B Larremore, Mirta Galesic, Jennie~E
  Brand, and Aaron Clauset.
\newblock Socioeconomic roots of academic faculty.
\newblock {\em Nature Human Behaviour}, pages 1--9, 2022.

\bibitem{myers2011mathematical}
Sean~A Myers, Peter~J Mucha, and Mason~A Porter.
\newblock Mathematical genealogy and department prestige.
\newblock {\em Chaos: An Interdisciplinary Journal of Nonlinear Science},
  21(4):041104, 2011.

\bibitem{Mathfest2023}
Mathematical~Association of~America.
\newblock {MAA MathFest 2023 in Tampa, FL}.
\newblock https://www.maa.org/meetings/mathfest.

\bibitem{GuidelinesUM}
Mathematical~Association of~America.
\newblock Guidelines for programs and departments in undergraduate mathematical
  sciences, 2003.
\newblock
  https://www.maa.org/sites/default/files/pdf/guidelines/Dept-Guidelines-Feb2003.pdf.

\bibitem{SelfStudy}
Mathematical~Association of~America.
\newblock Guidelines for undertaking a self-study in the mathematical sciences,
  2010.
\newblock
  https://www.maa.org/sites/default/files/pdf/ProgramReview/MAA-SelfStudyManual.pdf.

\bibitem{CUPM}
Mathematical~Association of~America.
\newblock Curriculum guide to majors in the mathematical sciences, 2015.
\newblock
  https://www.maa.org/sites/default/files/pdf/ProgramReview/MAA-SelfStudyManual.pdf.

\bibitem{reys2022some}
Robert Reys, Barbara Reys, and Jeffrey Shih.
\newblock Some {P}atterns of {P}hds in mathematics {A}warded {A}nnually by
  {I}nstitutions of {H}igher {E}ducation in the {U}nited {S}tates over the
  {L}ast {T}wo {D}ecades.
\newblock {\em Notices of the American Mathematical Society}, 69(1), 2022.

\bibitem{rissler2020gender}
Leslie~J Rissler, Katherine~L Hale, Nina~R Joffe, and Nicholas~M Caruso.
\newblock Gender differences in grant submissions across science and
  engineering fields at the {NSF}.
\newblock {\em Bioscience}, 70(9):814--820, 2020.

\bibitem{roca2023disruptjmm}
Rachel Roca, Carrie~Diaz Eaton, Drew Lewis, Joseph~E Hibdon~Jr, and Stefanie
  Marshall.
\newblock {\#DisruptJMM}: Themes of online social justice advocacy and
  community building in {STEM}.
\newblock {\em Journal of Humanistic Mathematics}, 2023.

\bibitem{SaxeBraddy}
Karen Saxe and Linda Braddy.
\newblock A common vision for undergraduate mathematical sciences programs in
  2025, 2015.
\newblock https://www.maa.org/sites/default/files/pdf/CommonVisionFinal.pdf.

\bibitem{schlenker2020prestige}
Jean-Marc Schlenker.
\newblock The prestige and status of research fields within mathematics.
\newblock {\em arXiv preprint arXiv:2008.13244}, 2020.

\bibitem{schmaling2023gender}
Karen~B Schmaling and Stephen~A Gallo.
\newblock Gender differences in peer reviewed grant applications, awards, and
  amounts: a systematic review and meta-analysis.
\newblock {\em Research Integrity and Peer Review}, 8(1):2, 2023.

\bibitem{storage2016frequency}
Daniel Storage, Zachary Horne, Andrei Cimpian, and Sarah-Jane Leslie.
\newblock The frequency of “brilliant” and “genius” in teaching
  evaluations predicts the representation of women and african americans across
  fields.
\newblock {\em PloS one}, 11(3):e0150194, 2016.

\bibitem{taffe2021racial}
Michael~A Taffe and Nicholas~W Gilpin.
\newblock Racial inequity in grant funding from the us national institutes of
  health.
\newblock {\em Elife}, 10:e65697, 2021.

\bibitem{topaz2020comparing}
Chad~M Topaz, James Cart, Carrie Diaz~Eaton, Anelise Hanson~Shrout, Jude~A
  Higdon, Kenan Ince, Brian Katz, Drew Lewis, Jessica Libertini, and
  Christian~Michael Smith.
\newblock Comparing demographics of signatories to public letters on diversity
  in the mathematical sciences.
\newblock {\em PloS one}, 15(4):e0232075, 2020.

\bibitem{topaz2020institute}
Chad~M Topaz, Maria-Veronica Ciocanel, Phoebe Cohen, Miles Ott, and Nancy
  Rodr{\'\i}guez.
\newblock Institute for the quantitative study of inclusion, diversity, and
  equity (qside).
\newblock {\em Notices of the American Mathematical Society}, 67(2), 2020.

\bibitem{topaz2016gender}
Chad~M Topaz and Shilad Sen.
\newblock Gender representation on journal editorial boards in the mathematical
  sciences.
\newblock {\em PLoS One}, 11(8):e0161357, 2016.

\bibitem{Wellcome}
The~Wellcome Trust.
\newblock Our commitment to tackling racism at wellcome, 2020.
\newblock
  https://wellcome.org/press-release/our-commitment-tackling-racism-wellcome.

\bibitem{wang2021science}
D.~Wang and A.L. Barab{\'a}si.
\newblock {\em The {S}cience of {S}cience}.
\newblock Cambridge University Press, 2021.

\bibitem{wapman2022quantifying}
K~Hunter Wapman, Sam Zhang, Aaron Clauset, and Daniel~B Larremore.
\newblock Quantifying hierarchy and dynamics in {US} faculty hiring and
  retention.
\newblock {\em Nature}, 610(7930):120--127, 2022.

\bibitem{way2016gender}
Samuel~F Way, Daniel~B Larremore, and Aaron Clauset.
\newblock Gender, productivity, and prestige in computer science faculty hiring
  networks.
\newblock In {\em Proceedings of the 25th international conference on world
  wide web}, pages 1169--1179, 2016.

\bibitem{zeng2017science}
An~Zeng, Zhesi Shen, Jianlin Zhou, Jinshan Wu, Ying Fan, Yougui Wang, and
  H~Eugene Stanley.
\newblock The science of science: From the perspective of complex systems.
\newblock {\em Physics reports}, 714:1--73, 2017.

\end{thebibliography}

%%%%%%%%%%%%%%%%%%%%%%%%%%%%%%%%%%%%
\section*{Authors}
\begin{itemize}
    \item \emph{Ron Buckmire} received his Ph.D. in mathematics from Rensselaer Polytechnic Institute in 1994 and has been on the faculty of Occidental College in Los Angeles, California ever since. He is a passionate advocate for broadening the participation of people from historically marginalized and currently under-represented groups in the mathematical sciences. 
Dr. Buckmire believes that: 1) mathematics is a human endeavor; 2) mathematics is created, discovered, learned, researched, and taught by people; 3) the identities and experiences of ``who does math" are important and that everyone and anyone can and should be welcome in the mathemaitcs community.

Mathematics Department, Occidental College, Los Angeles, CA \\
ron@oxy.edu

\item \emph{Joseph E. Hibdon, Jr.} received his PhD in Applied Mathematics from Northwestern University in 2011.  He is an associate professor of mathematics at Northeastern Illinois University, a Hispanic Serving Institutions in Chicago, Illinois.  Dr. Hibdon works in many interdisciplinary teams across the country to increase diversity in STEM and broadening the participation of students in research at the undergraduate level.  Dr. Hibdon's recent research is in mathematical modeling, with a focus on public health and biological phenomenon.

Department of Mathematics, Northeastern Illinois University, Chicago, Illinois\\
j-hibdonjr@neiu.edu

\item \emph{Drew Lewis} received his Ph.D. in mathematics from Washington University in St. Louis in 2012 and most recently served on the faculty at the University of South Alabama. He now works primarily in education research and faculty development.

drew.lewis@gmail.com

\item \emph{Omayra Ortega} received her PhD in Applied Mathematics \& Computational Sciences from the University of Iowa in 2008. She is an Associate Professor of Applied Mathematics \& Statistics and the Assistant Dean of Research and Internships in the School of Science and Technology at Sonoma State University. Using the tools from statistics, mathematics, data science, public health and epidemiology, Dr. Ortega tackles the emerging health issues. She is deeply committed to broadening the participation of under-represented minorities in STEM and mentoring students through the challenges of academia.

Department of Mathematics \& Statistics, Sonoma State University, Rohnert Park, CA\\ ortegao@sonoma.edu

\item \emph{Jos\'e L. Pab\'on} is a Ph.D. candidate in mathematics at the New Jersey Institute of Technology in Newark, N.J. He received his A.B. degree in mathematics from Princeton University in 2019 and has a particular affinity for working with disenfranchised populations, having first-hand experience with these given he was born and raised in Puerto Rico.

Mathematical Sciences Department, New Jersey Institute of Technology, Newark, NJ\\
jlp43@njit.edu

\item \emph{Rachel Roca} is a Ph.D. student in computational mathematics, science, and engineering at Michigan State University. Her interests lie in topological data analysis and computing education with an emphasis in social good. She believes in leveraging mathematical and data science tools for both advocacy and activism.

Department of Computational Mathematics, Science, and Engineering, Michigan State University, East Lansing, MI \\
rocarach@msu.edu

\item \emph{Andr\'es R. Vindas-Mel\'endez} received their PhD in mathematics from the University of Kentucky in 2021. 
Andr\'es is an NSF Postdoctoral Fellow and Lecturer at UC Berkeley and will begin as a tenure-track Assistant Professor at Harvey Mudd College in July 2024. 
Andr\'es' research interests are in algebraic, enumerative, and geometric combinatorics and have expanded to include applications of data science and mathematics for social justice. 

Department of Mathematics, University of California, Berkeley, CA \\
Department of Mathematics, Harvey Mudd College, Claremont, CA \\
andres.vindas@berkeley.edu; arvm@hmc.edu

\end{itemize}

%%%%%%%%%%%%%%%%%%%%%%%%%%%%%%%%%%%%

\end{document}